\documentclass[a4paper,11pt]{article}
\usepackage{amssymb}

\newcommand{\R}{\mathbb{R}}

\newcommand{\N}{\mathbb{N}}

\newcommand{\beq}{\begin{equation} }
\newcommand{\eqq}{\end{equation} }
\newcommand{\cuad}{{\sqcap\kern-.68em\sqcup}}

\newcommand{\norm}[1]{\|#1\|}

\newtheorem{definition}{Definition}[section]
\newtheorem{teo}{Theorem}[section]

\newtheorem{proposition}{Proposition}[section]

\newtheorem{lemma}{Lemma}[section]

\newtheorem{remark}{Remark}[section]
\newcommand{\bremark}{\begin{remark} \em}
\newcommand{\eremark}{\end{remark} }

\def\beeq{\begin{equation}}
\def\eeq{\end{equation}}
\newcommand{\begeqaet}{\begin{eqnarray*}}
\newcommand{\eneqaet}{\end{eqnarray*}}

\begin{document}
\begin{center}{\bf  \Large
On a class of semilinear fractional elliptic equations\medskip

 involving outside Dirac data }

  \bigskip  \medskip

 Huyuan Chen\footnote{hc64@nyu.edu}  \qquad  Hichem Hajaiej\footnote{hh62@nyu.edu}
  \qquad  Ying Wang\footnote{yingwang00@126.com}
\medskip

\begin{abstract}
The purpose of this article is to
give a complete study of
the  weak solutions of the fractional elliptic equation
\begin{equation}\label{00}
 \arraycolsep=1pt
\begin{array}{lll}
 (-\Delta)^{\alpha} u+u^p=0\ \ \ \ &\ {\rm in}\ \ B_1(e_N),\\[2mm]\phantom{(-\Delta)^{\alpha} +u^p}
 u=\delta_{0}& \ {\rm in}\ \ \R^N\setminus B_1(e_N),
\end{array}
\end{equation}
where  $p\ge0$, $ (-\Delta)^{\alpha}$ with $\alpha\in(0,1)$  denotes  the fractional Laplacian operator in the principle value sense,  $B_1(e_N)$ is the unit ball
 centered at $e_N=(0,\cdots,0,1)$ in $\mathbb{R}^N$ with $N\ge 2$ and $\delta_0$ is the Dirac mass concentrated at the origin. We prove that  problem (\ref{00}) admits a unique weak solution when
 $p> 1+\frac{2\alpha}{N}$. Moreover, if in addition $p\ge \frac{N+2}{N-2}$, the weak solution vanishes as $\alpha\to 1^-$. We also show that
problem (\ref{00}) doesn't have any weak solution  when $p\in[0, 1+\frac{2\alpha}{N}]$.
These results are very surprising since there are in total contradiction with the classical setting, i.e.
$$
\arraycolsep=1pt
\begin{array}{lll}
 -\Delta u+ u^p=0\ \ \ \ &\ {\rm in}\  \  B_1(e_N),\\[2mm]
 \phantom{-\Delta  +u^{p} }
u=\delta_{0}& \ {\rm in}\ \ \R^N\setminus B_1(e_N),
\end{array}
$$
for which it has been proved that there are no solutions for $p\ge \frac{N+1}{N-1}$.

\end{abstract}
\end{center}

\vspace{1mm}

  \smallskip



\setcounter{equation}{0}
\section{ Introduction}
Fractional PDEs have gained tremendous interest, not only from mathematicians but also from physicists and engineering,
during the last years. This is essentially due to their widespread domains of applications. In fact the fractional Laplacian arises in many
areas including medicine \cite{E}, bio-engineering \cite{M1,M2,M3,M4}, relativistic physics\cite{BDS,H,LY}, Modeling populations \cite{TZ}, flood flow, material viscoelastic theory, biology and earthquakes.   It is also particularly relevant to study some situations, in which  the fractional Laplacian is involved in PDEs, featuring  irregular data such that those phenomena describing source terms which are concentrated at points.  In our context, the source is placed outside the unit ball $B_1(e_N)$. This generates long-term interactions and short-term interactions, described by the nonlocal operator $(-\Delta)^\alpha$  and the nonlinear absorption $u^p$ respectively. $(-\Delta)^\alpha$ has also a probabilistic interpretation, related to the above one. It is the $\alpha-$stable subordinated infinitesimal killed Brownian motion.

Let $B_1(e_N)$ be the unit ball in
$\mathbb{R}^N(N\ge2)$ with  center $e_N=(0,\cdots,0,1)$ and $\delta_{0}$ be the Dirac mass
concentrated at the origin. Our main objective in this article is to investigate the existence, nonexistence and uniqueness of positive weak solutions of the
semilinear fractional equation
\begin{equation}\label{1.1}
 \arraycolsep=1pt
\begin{array}{lll}
 (-\Delta)^{\alpha} u+ u^p=0\ \ \ \ &\ {\rm in}\ \
B_1(e_N),\\[2mm]\phantom{(-\Delta)^{\alpha} +u^{p} }
u=\delta_{0}& \ {\rm in}\ \ \R^N\setminus B_1(e_N),
\end{array}
\end{equation}
where $p\ge0$  and the fractional
Laplacian $ (-\Delta)^{\alpha}$ with $\alpha\in(0,1)$ is defined by
$$ (-\Delta)^\alpha  u(x)=c_{N,\alpha}\lim_{\epsilon\to0^+} (-\Delta)_\epsilon^\alpha u(x),$$
where \begin{equation}\label{1.2}
  c_{N,\alpha}=\left(\int_{\R^N} \frac{1-\cos(z_1)}{|z|^{N+2\alpha}}dz\right)^{-1}
\end{equation}
 with $z=(z_1,\cdots,z_N)\in \R^N$ and
$$
(-\Delta)_\epsilon^\alpha  u(x)=-\int_{\R^N\setminus B_\epsilon(x)}\frac{ u(z)-
u(x)}{|z-x|^{N+2\alpha}} dz.
$$

In 1991, a fundamental contribution to semilinear elliptic
equations involving measures as boundary data is due to Gmira and V\'{e}ron  \cite{GV}, where they  studied the existence and uniqueness of weak solutions for
\begin{equation}\label{1.1.1}
\arraycolsep=1pt
\begin{array}{lll}
-\Delta  u+h(u)=0\quad &{\rm in}\quad \Omega,
\\[2mm]\phantom{-----}
u=\mu\quad&{\rm on}\quad \partial\Omega,
\end{array}
\end{equation}
where $\Omega$ is a bounded  $C^2$ domain and $\mu$ is a bounded Radon measure defined in $\partial\Omega$.
A function $u$ is said to be a weak solution of (\ref{1.1.1})
{\it  if $u\in
L^1(\Omega)$, $h(u)\in L^1(\Omega,\rho dx)$  and
\begin{equation}\label{1.1.1.0}
\int_\Omega [u(-\Delta)\xi+ h(u)\xi]dx=\int_{\partial\Omega}\frac{\partial\xi(x)}{\partial\vec{n}_x}d\mu(x),\quad \forall\xi\in C^{1.1}_0(\Omega),
\end{equation}
where $\rho(x)={\rm dist}(x,\partial\Omega)$ and $\vec{n}_x$ denotes the unit inward normal vector  at a point $x$.}
Gmira and V\'{e}ron  proved that the problem (\ref{1.1.1}) admits a unique weak solution  when
 $h$ is a continuous and nondecreasing function satisfying
\begin{equation}\label{14.04}
\int_1^\infty [h(s)-h(-s)]s^{-1-\frac{N+1}{N-1}}ds<+\infty.
\end{equation}
The weak solution of (\ref{1.1.1}) is approached by the classical solutions  of (\ref{1.1.1}) when $\mu$ is replaced by a
sequence of regular functions $\{\mu_n\}$, which converge to $\mu$ in the distribution sense. Furthermore, they showed that there is no weak solution of (\ref{1.1.1})
  when  $\mu=\delta_{x_0}$ with $x_0\in\partial\Omega$ and $h(s)=|s|^{p-1}s$ with
$p\ge\frac{N+1}{N-1}$. Later on, this subject has been vastly expanded in recent works,
see the papers of Marcus and V\'{e}ron \cite{MV1,MV2,MV3,MV4}, Bidaut-V\'{e}ron  and  Vivier \cite{BY} and references therein.

In the fractional setting, the equivalent of  (\ref{1.1.1})  when  $\mu=\delta_{x_0}$ has been considered in \cite{CH}, where the authors proved that the
weak solution of
\begin{equation}\label{eq 1.0}
\arraycolsep=1pt
\begin{array}{lll}
 (-\Delta)^\alpha   u+ g(u)=k\frac{\partial^\alpha\delta_{x_0}}{\partial \vec{n}^\alpha}\quad  &{\rm in}\quad\ \ \bar\Omega,\\[3mm]
 \phantom{------\ \  }
 u=0\quad &{\rm in}\quad\ \ \bar\Omega^c
 \end{array}
\end{equation}
is approximated by the weak solutions, as $t\to0^+$, of
$$
\arraycolsep=1pt
\begin{array}{lll}
 (-\Delta)^\alpha   u+ g(u)=kt^{-\alpha}\delta_{x_0+t\vec{n}_{x_0}}\quad  &{\rm in}\quad\ \ \bar\Omega,\\[3mm]
 \phantom{------\ \  }
 u=0\quad &{\rm in}\quad\ \ \bar\Omega^c.
 \end{array}
$$
More precisely, in the fractional setting, $\frac{\partial^\alpha\delta_{x_0}}{\partial \vec{n}^\alpha}$ plays the same role of
$u=\delta_{x_0}$ on $\partial\Omega$ in (\ref{1.1.1}). Our purpose in this article is to study the solution of
$$(-\Delta)^{\alpha} u+ u^p=0\quad\ {\rm in}\quad B_1(e_N)$$
when  the exact Dirac mass concentrated at the origin is considered. Our main idea is to  make use of nonlocal properties of the fractional Laplacian
to move the Dirac mass at $-te_N$ when $t\to0^+$ and we then proceed by approximation techniques. Before giving our main results, we must first give an appropriate definition of  weak solution of (\ref{1.1}).  It is then worth to mention two important results. The equation
 \begin{equation}\label{eq 1.1}
  \arraycolsep=1pt
\begin{array}{lll}
 (-\Delta)^{\alpha} u+u^p=0\ \ \ \ &\ {\rm in}\ \
B_1(e_N),\\[2mm]\phantom{-----\ }
u=f& \ {\rm in}\ \ \R^N\setminus B_1(e_N),
\end{array}
\end{equation}
where $f\in C_0(\R^N\setminus B_1(e_N))$,
admits a unique classical solution $u_f$, see \cite[Theorem 2.5]{CFQ}. Furthermore, let $\tilde u_f=u_f$ in $B_1(e_N)$ and $\tilde u_f=0$ in $\R^N\setminus B_1(e_N)$,
then $\tilde u_f$ is the unique classical solution of
 \begin{equation}\label{eq 1.2}
 \arraycolsep=1pt
\begin{array}{lll}
 (-\Delta)^{\alpha} u+u^p=c_{N,\alpha}\int_{\R^N\setminus B_1(e_N)}\frac{f(y)}{|x-y|^{N+2\alpha}}dy\ \ \ \ &\ {\rm in}\ \
B_1(e_N),\\[2mm]\phantom{-----\ }
u=0& \ {\rm in}\ \ \R^N\setminus B_1(e_N)
\end{array}
\end{equation}
and satisfies the identity:
$$
\int_{B_1(e_N)}\left[u(x)(-\Delta)^\alpha \xi(x) +u^p(x)\xi(x)\right]dx = c_{N,\alpha}
\int_{B_1(e_N)}\int_{\R^N\setminus B_1(e_N)}\frac{\xi(x)f(y)}{|x-y|^{N+2\alpha}}dydx,
$$
for any $\xi\in C^\infty_0(B_1(e_N))$. Let us mention that $C_0^\infty(B_1(e_N))$ is the space of test functions
$\xi\in C^\infty(\R^N)$ with support in  $B_1(e_N)$.

Inspired by above identity, we give the definition of   weak solution to (\ref{1.1}) as follows.
\begin{definition}\label{def 3}
  We say that $u$ is a weak solution of (\ref{1.1}) if $u\in  L^1(B_1(e_N))$, $u^p\in L^1_{\rm loc}(B_1(e_N))$  and
\begin{equation}\label{1.7}
\int_{B_1(e_N)} \left[u(x)(-\Delta)^\alpha \xi(x) +u^p(x)\xi(x)\right]dx =
\int_{B_1(e_N)}\xi(x)\Gamma_0(x) dx,\quad \forall \xi\in C^\infty_0(B_1(e_N)),
\end{equation}
\end{definition}
where
\begin{equation}\label{Gamma}
  \Gamma_0(x)=\frac{c_{N,\alpha}}{|x|^{N+2\alpha}},\quad\forall x\in\R^N\setminus\{0\}.
\end{equation}

It is well known that the definition of the weak solution heavily depends on the test functions space and the best function space is the one
 which enables us to get the "strongest" weak solution. 
  In \cite{CV1,CV2}, semilinear fractional equations with measures has been studied via  the test functions  space $\mathbb{X}_{\alpha, \Omega}\subset C(\R^N)$ for a $C^2$ bounded open domain $\Omega$,   where {\it $\mathbb{X}_{\alpha, \Omega}$ is the space of functions
$\xi$ satisfying:\smallskip

 (1) ${\rm supp}(\xi)\subset\bar\Omega$; \smallskip

 (2) $(-\Delta)^\alpha\xi(x)$ exists for all $x\in \Omega$
and $|(-\Delta)^\alpha\xi(x)|\leq C$ for some $C>0$;\smallskip

 (3) there exist $\varphi\in L^1(\Omega,\rho^\alpha dx)$
and $\epsilon_0>0$ such that $|(-\Delta)_\epsilon^\alpha\xi|\le
\varphi$ a.e. in $\Omega$, for all
$\epsilon\in(0,\epsilon_0]$,

where $\rho(x)=dist(x,\partial\Omega)$.
}\\[0.5mm]
The test functions space $C_0^\infty(B_1(e_N))$ has  stronger topology than $\mathbb{X}_{\alpha,B_1(e_N)}$ does,
 the weak solution in Definition 1.1 with test functions space $\mathbb{X}_{\alpha,B_1(e_N)}$ would be
stronger than the one with test functions space $C_0^\infty(B_1(e_N))$. It is then worth
to mention that the test functions space $C_0^\infty(B_1(e_N))$ could not be replaced to
the test functions space $\mathbb{X}_{\alpha,B_1(e_N)}$ in our setting. For example,  $\xi_0:=\mathbb{G}_\alpha[1]\in \mathbb{X}_{\alpha,B_1(e_N)}$, but
(\ref{1.7}) does not hold for  $\xi_0$,
where $G_\alpha$ denotes the Green kernel of $(-\Delta)^\alpha$ in $B_1(e_N)\times B_1(e_N)$ and  $\mathbb{G}_\alpha$ is the Green
operator defined as
\begin{equation}\label{1.8}
 \mathbb{G}_\alpha[f](x)=\int_{B_1(e_N)} G_\alpha(x,y)f(y)dy,\quad f\in L^1(B_1(e_N),\rho^\alpha dx).
\end{equation}

Let us state our existence result.

\begin{teo}\label{teo 2}
Assume that $\alpha\in(0,1)$ and $p>1+\frac{2\alpha}N$.
Then  there exists   a unique nonnegative weak solution $u_{\alpha,p}$ of
 (\ref{1.1}) such that for some $c_1>1$, we have
\begin{equation}\label{1.10}
 0<u_{\alpha,p}(x)\le c_1   |x|^{-\frac{N+2\alpha}p},\qquad\forall x\in B_1(e_N)
\end{equation}
and
\begin{equation}\label{1.11}
\frac{1}{c_1}t^{-\frac{N+2\alpha}p}\le u_{\alpha,p}(te_N) \le
c_1 t^{-\frac{N+2\alpha}p},\qquad\forall t\in (0,1).
\end{equation}

\end{teo}

\begin{remark}\label{remark 29111}
$(i)$ The existence results are very surprising as they are in total different from the Laplacian case,
where   (\ref{1.1.1}) with $\mu=\delta_0$ has a weak solution only when $p<\frac{N+1}{N-1}$.

$(ii)$ From (\ref{1.11}), the singularity is only near the origin.  We also notice that
$$u_{\alpha,p}^p(te_N)\ge \frac1{c_1} t^{-(N+2\alpha)},\qquad\forall t\in(0,1),$$
which implies that the absorption nonlinearity $u^p$ plays a primary role in  the equation (\ref{1.1}).
While the absorption nonlinearity always plays a second role in a measure  framework.

$(iii)$ The uniqueness cannot directly follow   Kato's inequality \cite[Proposition 2.4]{CV1} since it has been built
in the framework of the test functions space $\mathbb{X}_{\alpha,B_1(e_N)}$. In this paper, as mentioned above,
$C_0^\infty(B_1(e_N))$ is the appropriate test functions space. This will give birth to a lot of technical difficulties to prove
the existence, nonexistence and uniqueness of weak solutions of (\ref{1.1}).
\end{remark}

If $p>1+\frac{2\alpha}N$, the weak solution $u_{\alpha, p}$ of (\ref{1.1}) is approximated by the unique solution
$u_s$  $(s\in(0,1))$
 of
\begin{equation}\label{eq 1.3}
  \arraycolsep=1pt
\begin{array}{lll}
 (-\Delta)^{\alpha} u +u^p=0\quad   &{\rm in}
\quad  B_1(e_N),\\[2mm]\phantom{-----\ }
u=\delta_{-se_N}\quad &  {\rm in}\quad \R^N\setminus B_1(e_N).
\end{array}
\end{equation}

When  $p\in [0,1+\frac{2\alpha}N]$, we will prove that $\{u_s\}$ blows up everywhere in $B_1(e_N)$ as $s\to0^+$,
therefore,  we can deduce the  nonexistence of weak solutions of (\ref{1.1}) when $p\le1+\frac{2\alpha}N$.
More precisely, we have the following results.

\begin{teo}\label{teo 3}
Assume that $\alpha\in(0,1)$ and $0\leq p\le1+\frac{2\alpha}N$.
Then problem (\ref{1.1}) does not have any weak solution.
\end{teo}

In the proof of Theorem \ref{teo 3}, we will first need to prove the crucial estimate
$$u_{\alpha,p}^p(x)\ge c_2|x|^{-(N+2\alpha)},\quad\forall x\in \mathcal{C},$$
where $\mathcal{C}=\{x\in\R^N:\  \exists t\in(0,1)\ {\rm s.t.}\  |x-te_N|<\frac t8 \}$ is a cone in $B_1(e_N)$. We combine the symmetry property and decreasing monotonicity
in our proof of the nonexistence. This phenomena is due to the nonlocal characteristic of fractional Laplacian
that requires the functions to be in $L^1_{\rm loc}(\R^N)$.

Finally, our interest is to study the asymptotic behavior of $u_{\alpha,p}$ as $\alpha$ goes to $1^-$.
\begin{teo}\label{teo 4}
Assume that $\alpha\in(0,1)$, $p\ge\frac{N+1}{N-1}$ and $ u_{\alpha,p} $ is the unique weak solution of problem (\ref{1.1}).
Then   $\{u_{\alpha,p}\}_\alpha$  vanishes as $\alpha\to1^-$.
\end{teo}

\begin{remark}
 $(i)$ for $p\ge  \frac{N+2}{N-2}$, a sequence of barrier functions, which converge to 0 locally in $B_1(0)$, could be constructed directly
to control $\{u_{\alpha,p}\}_\alpha$;\\
$(ii)$ for $\frac{N+1}{N-1}\le p<\frac{N+2}{N-2}$,   Theorem 3.1 in \cite{GV} is involved to control  $\{u_{\alpha,p}\}_\alpha$;\\
$(iii)$ for $p<1+\frac{2}N$, there exists $\alpha_p\in(0,1)$ such that $p\le 1+\frac{2\alpha}N$ for $\alpha\in(\alpha_p,1)$, there is no weak solution for problem (\ref{1.1}) from Theorem \ref{teo 3}; \\
$(iv)$ for $p\in [1+\frac2N, \frac{N+1}{N-1})$,  it is still open for the limit of $\{u_{\alpha,p}\}_\alpha$ as $\alpha\to1^-$.
\end{remark}

In Section 2, we treat the problem (\ref{eq 1.3}). When the Dirac mass concentrates at point $-se_N$
away from $\bar \Omega$, we build the existence, uniqueness weak solution $u_s$ of (\ref{eq 1.3}) and show how the Dirac mass
is transformed into the nonhomogeneous term. In this case, the test functions space could be improved
into $\mathbb{X}_{\alpha,B_1(e_N)}$, since the solution has no singularity in $\bar\Omega$.

In Section 3, we give a detailed account of  the procedure enabling us to move the singular points   $\{-se_N\}$  to the origin. The first difficulty
arises from the fact that $\mathbb{G}_\alpha[\Gamma_s]$ blows up everywhere as $s\to0^+$ [see Lemma 3.1], that is,
there is no solution of
\begin{equation}\label{eq 1.4}
  \arraycolsep=1pt
\begin{array}{lll}
 (-\Delta)^{\alpha} u =0\quad   &{\rm in}
\quad  B_1(e_N),\\[2mm]\phantom{---\ }
u=\delta_{0}\quad &  {\rm in}\quad \R^N\setminus B_1(e_N).
\end{array}
\end{equation}
Therefore, we have to resort a barrier function, that is the minimal classical solution of
$$ \arraycolsep=1pt
\begin{array}{lll}
 (-\Delta)^{\alpha} u+u^p=\Gamma_0\ \ \ \ & {\rm in}\ \
B_1(e_N),\\[2mm]\phantom{(-\Delta)^{\alpha}+u^p}
u=0& {\rm in} \ \ B_1^c(e_N)\setminus\{0\}.
\end{array}
$$
In order to control the limit of $\{u_s\}$ near $\partial B_1$, especially near the origin, some typical  truncated functions
have to be constructed carefully.  The second difficulty comes from the proof of the uniqueness.  We proceed by  contradiction, assuming that there is two
solutions  and we will show their difference   could be   improved the test function from $C^\infty_0$ into $\mathbb{X}_{\alpha,B_1(e_N)}$,  this
enables us to use Kato's inequality \cite[Proposition 2.4]{CV1} and to conclude.

Section 4 is devoted to blow-up case. The difficulty is to obtain the blow-up everywhere in $B_1(e_N)$ just from a lower bounds of $u_s$, see Lemma 4.1.
To overcome it, we combine the symmetric of the domain and resort the symmetry result of $u_s$ and then one point blowing up leads
to blowing up every where in $B_1(e_N)$.

Finally, we analyse  decay approximation of the weak solution for problem (\ref{1.1}) when $\alpha\to1^-$. For $p\ge \frac{N+2}{N-2}$, the  first challenge  is to construct a sequence upper bounds that converges to zero. To this end, we have to study $\lim_{\alpha\to1^-}(-\Delta)^\alpha \Phi_\sigma$, where $\Phi_\sigma(x)=|x|^{-\sigma}$ and then use proper parameters to construct the bounds.  For $\frac{N+1}{N-1}\le p\le \frac{N+2}{N-2}$, $\Phi_\sigma$ could not be used to construct properly the upper bounds, and then we use some argument of \cite{GV}.

\setcounter{equation}{0}
\section{Dirac mass concentrated at $\{-se_N\}$ with $s\in(0,1)$ }







The purpose of this section is to introduce some preliminaries.
First let us state an important Comparison Principle.

\begin{teo}\cite[Theorem 2.3]{CFQ}\label{comparison} \ \
Let $u$ and $v$ be  super-solution
and sub-solution, respectively, of
$$
 (-\Delta)^{\alpha} u+h(u)=f\ \ {\rm in}\ \ \mathcal{O},$$ where $\mathcal{O}$ is an open, bounded and connected
domain of class $C^2$, the function $f:\mathcal{O}\to\R$ is
continuous and $h:\R\to\R$ is increasing.

Suppose that $v(x)\le u(x),\ \forall x\in\mathcal{O}^c$, $u$ and $v$ are
continuous in $\bar \mathcal{O}$. Then
$$u(x)\ge v(x),\quad \forall x\in \mathcal{O}.$$

\end{teo}

Now we investigate the weak solution of
\begin{equation}\label{eq 2.1}
  \arraycolsep=1pt
\begin{array}{lll}
 (-\Delta)^{\alpha} u +u^p=0\quad   &{\rm in}
\quad  B_1(e_N),\\[2mm]\phantom{-----\ }
u=\delta_{-se_N}\quad &  {\rm in}\quad \R^N\setminus B_1(e_N),
\end{array}
\end{equation}
where $s\in(0,1)$.  To this end, we construct a sequence of $C^2$ functions to approximate the
Dirac measure.  Let $g_0:\R^N\to[0,1]$ be a radially symmetric decreasing $C^2$ function  with the support in $\overline{B_{\frac12}(0)}$
such that $\int_{\R^N} g_0(x)dx=1$.
For any $n\in\N$ and $s\in (0,1)$, we denote
$$g_n(x)=n^{N}g_0(n(x+se_N)),\qquad\forall x\in \R^N.$$
Then we certainly have  that
$$
g_n\rightharpoonup\delta_{-se_N}\ \ {\rm as}\ \ n\to+\infty,
$$
in the distribution sense and for any $s>0$, there exists $N_s>0$ such that for any $n\ge N_s$,
$${\rm supp} (g_n)\subset \overline{B_\frac s2(-se_N)}.$$

 In order to investigate
the solution of (\ref{eq 2.1}), we consider the approximating solution
$w_n$ of
\begin{equation}\label{eq 2.2}
  \arraycolsep=1pt
\begin{array}{lll}
 (-\Delta)^{\alpha} u +u^p =0\quad   &{\rm in}
\quad  B_1(e_N),\\[2mm]\phantom{-----\ }
u=g_n\quad &  {\rm in}\quad \R^N\setminus B_1(e_N).
\end{array}
\end{equation}
\begin{lemma}\label{lm 2.1}
Assume that $p>0$ and $\{g_n\}$ is a sequence of $C^2$ functions converging to $\delta_{-se_N}$ with supports in $\overline{B_\frac s2(-se_N)}$.
Denote that
\begin{equation}\label{2.2}
\tilde g_n(x):=c_{N,\alpha} \int_{\R^N}\frac{g_n(y)}{|x-y|^{N+2\alpha}}dy,\quad \forall x\in B_1(e_N).
\end{equation}
Then problem (\ref{eq 2.2}) admits a unique solution $w_n$ such that
$$0<w_n \le \mathbb{G}_\alpha[\tilde g_n]\quad{\rm in}\quad B_1(e_N). $$
Moreover, the function $\tilde w_n:=w_n\chi_{B_1(e_N)}$ is the unique solution of
\begin{equation}\label{eq 2.3}
  \arraycolsep=1pt
\begin{array}{lll}
 (-\Delta)^{\alpha} u +u^p =\tilde g_n\quad   &{\rm in}
\quad  B_1(e_N),\\[2mm]\phantom{-----\ }
u=0\quad &  {\rm in}\quad \R^N\setminus B_1(e_N).
\end{array}
\end{equation}

\end{lemma}
{\bf Proof.} The existence and uniqueness of solution to problem (\ref{eq 2.2}) refers to \cite[Theorem 2.5]{CFQ}.
For $n\ge N_s$, we have that supp$(g_n)\subset B_{\frac s2}(-se_N)$ and then $\tilde g_n\in C^1(\overline{B_1(e_N)})$ and
  $$\tilde w_n=w_n-g_n\quad {\rm in}\quad \R^N.$$
By the definition of fractional Laplacian, it implies that
\begin{eqnarray*}
  (-\Delta)^\alpha \tilde w_n(x)+\tilde w_n(x)^p&=&(-\Delta)^\alpha w_n(x)-(-\Delta)^\alpha g_n(x)+w_n(x)^p \\
   &=& c_{N,\alpha}\int_{\R^N}\frac{g_n(z)}{|z-x|^{N+2\alpha}} dz=\tilde g_n(x).
\end{eqnarray*}
Therefore, $\tilde w_n$ is a classical solution of (\ref{eq 2.3}) and
$$\tilde w_n\le \mathbb{G_\alpha}[\tilde g_n] \quad {\rm in}\quad B_1(e_N),$$
which implies that
$$w_n \le \mathbb{G}_\alpha[\tilde g_n]\quad{\rm in} \quad B_1(e_N).$$
The proof ends.\qquad$\Box$
\medskip

We remark that $\tilde w_n$ is the classical solution of (\ref{eq 2.3}), then by Lemma 2.1 and Lemma 2.2 in \cite{CV1}, we have that
\begin{equation}\label{2.4}
\int_{B_1(e_N)}[w_n(-\Delta)^\alpha \xi+w_n^p\xi]dx=\int_{B_1(e_N)} \xi\tilde g_n dx,\qquad \forall\xi\in C_0^\infty(B_1(e_N)).
\end{equation}
Here (\ref{2.4}) holds even for $\xi\in \mathbb{X}_{\alpha,B_1(e_N)}$.

\begin{lemma}\label{lm 2.2}
Let $\{\tilde g_n\}$ be defined in (\ref{2.2}), then $\tilde g_n$ converges to $\Gamma_s$
uniformly in $B_1(e_N)$ and in $C^{\theta}(B_1(e_N))$ for $\theta\in(0,1)$,
where
\begin{equation}\label{22-09-3}
\Gamma_s(x)=\frac{c_{N,\alpha}}{|x+se_N|^{N+2\alpha}},\quad \forall x\in\R^N\setminus \{-se_N\}.
\end{equation}

\end{lemma}
{\bf Proof.}   It is obvious that $\tilde g_n$ converges to $\Gamma_s$ every point in $ \overline{B_1(e_N)}$.
For $x,y\in B_1(e_N)$ and any $n\in\N$, we have that
\begin{eqnarray*}
 |\tilde g_n(x)-\tilde g_n(y)| &=&c_{N,\alpha}|\int_{B_{\frac s2}(-se_N)}[\frac{1}{|x-z|^{N+2\alpha}}-\frac{1}{|y-z|^{N+2\alpha}}]g_n(z)dz|
 \\ &\le & c_{N,\alpha}\int_{B_{\frac s2}(-se_N)}\frac{||x-z|^{N+2\alpha}- |y-z|^{N+2\alpha}|}{|x-z|^{N+2\alpha}|y-z|^{N+2\alpha}} g_n(z)dz
 \\ &\le & c_{N,\alpha}(N+2\alpha)|x-y|\int_{B_{\frac s2}(-se_N)}\frac{|x-z|^{N+2\alpha-1}+ |y-z|^{N+2\alpha-1}}{|x-z|^{N+2\alpha}|y-z|^{N+2\alpha}} g_n(z)dz
 \\ &\le & c_3|x-y|\int_{B_{\frac s2}(-se_N)}g_n(z)dz
 \\ &= & c_3|x-y|,
\end{eqnarray*}
where $c_3>0$ independent of $n$. So $\{\tilde g_n\}_n$ is uniformly bounded in $C^{0,1}(B_1(e_N))$. Combining the converging
$$\tilde g_n\to \Gamma_s\ {\rm every\ point\ in} \  \overline{B_1(e_N)}.$$
We conclude that $\tilde g_n$ converges to $\Gamma_s$
uniformly in $B_1(e_N)$ and in $C^{\theta}(B_1(e_N))$ for $\theta\in(0,1)$.
\qquad$\Box$

\begin{proposition}\label{pr 2.1}
Assume that $p>0$, $s\in (0,1)$ and $\Gamma_s$ is given by (\ref{22-09-3}). Then problem (\ref{eq 2.1}) admits a
unique weak solution $u_s$ such that
\begin{equation}\label{2.7}
 0\le u_s(x)\le \mathbb{G}_\alpha[\Gamma_s],\qquad x\in B_1(e_N).
\end{equation}
Moreover,  $\tilde u_s:=u_s\chi_{B_1(e_N)}$ is the unique classical solution of
\begin{equation}\label{eq 2.4}
 \arraycolsep=1pt
\begin{array}{lll}
 (-\Delta)^{\alpha} u +u^p =\Gamma_s\quad   &{\rm in}
\quad  B_1(e_N),\\[2mm]\phantom{-----\ }
u=0\quad &  {\rm in}\quad \R^N\setminus B_1(e_N).
\end{array}
\end{equation}

\end{proposition}
{\bf Proof. }  {\it Existence. }
It infers by Lemma \ref{lm 2.1}  that the solution $w_n$ of (\ref{eq 2.2}) satisfies that
\begin{equation}\label{2.6}
0<w_n \le \mathbb{G}_\alpha[\tilde g_n]\quad{\rm in}\quad B_1(e_N).
\end{equation}
By Lemma \ref{lm 2.2} we have that $\tilde g_n$ converges to $\Gamma_s$
uniformly in $\overline{B_1(e_N)}$ and in $C^\theta(B_1(e_N))$ with $\theta\in(0,1)$.
Therefore, there exists some constant $c_4>0$ independent of $n$ such that
$$\mathbb{G}_\alpha[\tilde g_n](x)\le \frac{c_4c_{N,\alpha}}{|x+se_N|^{N+2\alpha}}\le c_4c_{N,\alpha} s^{-N-2\alpha},\qquad\forall x\in B_1(e_N).$$
Thus,
$$\norm{w_n}_{L^\infty(B_1(e_N))}\le c_4c_{N,\alpha} s^{-N-2\alpha} ,\qquad  \norm{w_n}_{L^1(B_1(e_N))}\le c_4c_{N,\alpha} s^{-N-2\alpha}|B_1(e_N)|.$$

 By \cite[Theorem 1.2]{RS},  we have that
 \begin{equation}\label{1.3}
\begin{array}{lll}
\norm{\frac{w_n}{\rho^\alpha}}_{C^{\alpha}(  \overline{B_1(e_N)})} \le c_5[ \norm{w_n^p}_{L^{\infty}(B_1(e_N))}+\norm{\Gamma_s}_{L^{\infty}(B_1(e_N))}]\\[3mm]
\phantom{\norm{\frac{w_n}{\rho^\alpha}}_{C^{\alpha}(\bar B_1(e_N)}}
 \le c_6 [s^{-N-2\alpha}+s^{-(N+2\alpha)p}]
\end{array}
 \end{equation}
for some $c_5,c_6>0$.

In order to see the inner regularity, we  denote $\mathcal{O}_i$ the open sets with $i=1,2,3$ such that
 $$\mathcal{O}_1\subset \bar\mathcal{O}_1\subset \mathcal{O}_2\subset \bar\mathcal{O}_2\subset\mathcal{O}_3\subset \bar\mathcal{O}_3 \subset B_1(e_N). $$
 By  \cite[Lemma 3.1]{CV3}, for $\beta\in(0,\alpha)$,  there exists $c_7,c_8>0$ independent of $n$ such that
 $$\arraycolsep=1pt
\begin{array}{lll}
\norm{w_n}_{C^{\beta}(\mathcal{O}_2)} \le c_7[\norm{w_n}_{L^1(B_1(e_N))}+\norm{w_n^p}_{L^{\infty}(\mathcal{O}_3)}+\norm{w_n}_{L^{\infty}(\mathcal{O}_3)}]\\[3mm]
\phantom{\norm{w_n}_{C^{\beta}(\mathcal{O}_2)}}
 \le c_8 [s^{-N-2\alpha}+s^{-(N+2\alpha)p}].
\end{array}
$$
It follows by  \cite[Corollary 2.4]{RS} that there exist $c_9,c_{10}>0$ such that
\begin{equation}\label{2.0.10.0}
\arraycolsep=1pt
\begin{array}{lll}
\norm{w_n}_{C^{2\alpha+\beta}(\mathcal{O}_1)} \le {c_9}[\norm{w_n}_{L^1(B_1(e_N))}+\norm{w_n^p}_{C^{\beta}(\mathcal{O}_2)}
+\norm{w_n}_{C^{\beta}(\mathcal{O}_2)}]\\[3mm]
\phantom{\norm{w_n}_{C^{2\alpha+\beta}(\mathcal{O}_1)} }
 \le c_{10} [s^{-N-2\alpha}+s^{-(N+2\alpha)p}].
\end{array}
\end{equation}

Therefore,  by the Arzela-Ascoli Theorem, there exist $u_s\in C^{2\alpha+\epsilon}_{\rm loc}$ in $B_1(e_N)$ for some $\epsilon\in(0,\beta)$ and a subsequence $\{w_{n_k}\}$ such that
\begin{equation}\label{01-10}
w_{n_k}\to u_s\quad  {\rm locally\ in}\quad C^{2\alpha+\epsilon} \quad{\rm as}\quad  n_k\to\infty.
\end{equation}
Passing the limit of (\ref{2.4}) with $\xi\in \mathbb{X}_{\alpha,B_1(e_N)}$ as $n_k\to\infty$, we have that
\begin{equation}\label{2.5}
\int_{B_1(e_N)}[u_s(-\Delta)^\alpha \xi+u_s^p\xi]dx= \int_{B_1(e_N)} \xi(x)\Gamma_s(x).
\end{equation}
Moreover, since $w_n\to u_s$ and $\tilde g_n\to \Gamma_s$ uniformly in $B_1(e_N)$ as $n\to\infty$,
then it infers that
$$0\le u_s\le \mathbb{G}_\alpha[\Gamma_s]\quad {\rm in}\quad B_1(e_N).$$

{\it Uniqueness.} Let $v_s$ be a weak solution of (\ref{eq 2.2}) and then $\varphi_s:=u_{s}-v_s$ is a weak solution to
$$
\arraycolsep=1pt
\begin{array}{lll}
(-\Delta)^\alpha \varphi_s +u_{s}^p- v_s^p=0\quad &{\rm in}\quad B_1(e_N),
\\[2mm]\phantom{(-\Delta)^\alpha   +u_{s}^p- v_s^p}
\varphi_s=0\quad &{\rm in}\quad \R^N\setminus B_1(e_N).
\end{array}
$$
By Kato's inequality \cite[Proposition 2.4]{CV1},
$$\int_{B_1(e_N)}|\varphi_s|(-\Delta)^\alpha\xi+\int_{B_1(e_N)}[ u_s^p-v_s^p]{\rm sign}(u_s-v_s)\xi dx=0. $$
Taking $\xi=\mathbb{G}_\alpha[1]$, we have that
$$\int_{B_1(e_N)}[u_s^p- v_s^p]{\rm sign}(u_s-v_s)\xi dx\ge 0\quad{\rm and}\quad \int_{B_1(e_N)}|\varphi_s|dx=0,$$
then $\varphi_s=0$ a.e. in $B_1(e_N)$. Then the uniqueness is proved.

Furthermore,  we see that $\tilde w_n=w_n-g_n$ is  the unique classical solution of
\begin{equation} \label{4.9}
 \arraycolsep=1pt
\begin{array}{lll}
 (-\Delta)^{\alpha} u(x)+u^p(x)=\tilde g_n(x),\ \ \ \ &
\forall x\in B_1(e_N),\\[2mm]\phantom{(-\Delta)^{\alpha}+u^p(x)}
u(x)=0,& \forall x\in B_1(e_N)^c
\end{array}
\end{equation}
and $\tilde w_n$ converges to $\tilde u_s$ uniformly in $B_1(e_N)$.
By Stability Theorem   \cite[Lemma 4.5]{CS2} and (\ref{01-10}), $u_s\chi_{B_1(e_N)}$
is the classical solution of (\ref{eq 2.4}). \qquad$\Box$

\setcounter{equation}{0}
\section{Proof of Theorem \ref{teo 2}. }

In this section, we prove  Theorem \ref{teo 2} by moving the points $\{-se_N\}$ to the origin. To this end,  we
need derive more properties for $u_s$, where $u_s$ is the unique weak solution of (\ref{eq 2.1}).
\begin{lemma}\label{lm 3.1-26}
Let $p>0$, $s\in (0,1)$ and $u_s$ be the unique weak solution of (\ref{eq 2.1}). Then
the mapping $s\mapsto u_s$ is decreasing, that is,
$$u_{s_1}\ge u_{s_2}\quad {\rm if}\ s_1\le s_2. $$

\end{lemma}
{\bf Proof.} By Proposition \ref{pr 2.1},   $\tilde u_s:=u_s\chi_{B_1(e_N)}$ is the unique classical solution of
(\ref{eq 2.4}).

We claim that the mapping: $s\mapsto \Gamma_s$ is decreasing.
For $x\in B_1(e_N)$ and $s_1\le s_2$, we observe that $|x+s_1e_N|\le |x+s_2e_N|$,
then $\Gamma_{s_1}(x)\ge \Gamma_{s_2}(x).$  The claim is proved.

Therefore, for $s_1\le s_2$,  $u_{s_1}$ and $u_{s_2}$ are super solution and solution of (\ref{eq 2.4}) replaced $\Gamma_s$ by $\Gamma_{s_2}$,
then it infers by the Comparison Principle that $u_{s_1}\ge u_{s_2}$ in $B_1(e_N)$.  \qquad$\Box$


\begin{lemma}\label{lm 3.1}
Let   $s\in (0,1)$ and denote
$$ \mathbb{G}_\alpha[\Gamma_s](x)=\int_{B_1(e_N)}G_\alpha(x,y)\Gamma_s(y)dy,\quad \forall x\in B_1(e_N). $$
 Then
 $$\lim_{s\to0^+}\mathbb{G}_\alpha[\Gamma_s](x)=+\infty,\qquad \forall x\in  B_1(e_N).$$

\end{lemma}
{\bf Proof.} Using \cite[Theorem 1.2]{CS}, it follows that
$$G_\alpha(x,y)\ge \min\left\{\frac{c_{11} }{|x-y|^{N-2\alpha}},\frac{c_{11} \rho^\alpha(x)\rho^\alpha(y)}{|x-y|^N}\right\},\quad x,y\in B_1(e_N),$$
where $c_{11}>0$ dependent of ${N,\alpha}$.
Now for $x\in B_1(e_N)$ and $y\in B_1(e_N) \cap B_{\frac{|x|}4}(0)$, we have that
$$G_\alpha(x,y)\ge \frac{c_{11} \rho^\alpha(x)\rho^\alpha(y)}{|x-y|^N} $$
and
\begin{eqnarray*}
 \mathbb{G}_\alpha[\Gamma_s](x) &\ge & \int_{B_1(e_N)}\min\left\{\frac{c_{11}}{|x-y|^{N-2\alpha}},\frac{c_{11} \rho^\alpha(x)\rho^\alpha(y)}{|x-y|^N}\right\}\frac{c_{N,\alpha}}{|y+se_N|^{N+2\alpha}}dy \\&\ge &  \int_{B_1(e_N) \cap B_{\frac{|x|}4}(0)} \frac{c_{11} \rho^\alpha(x)\rho^\alpha(y)}{|x-y|^N}\frac{c_{N,\alpha}}{|y+se_N|^{N+2\alpha}}dy
 \\&\ge &\frac45 c_{11}\rho^\alpha(x)|x|^{-N}\int_{B_1(e_N) \cap B_{\frac{|x|}4}(0)} \frac{c_{N,\alpha} \rho^\alpha(y)}{|y+se_N|^{N+2\alpha}}dy
 \\&\to&+\infty \qquad{\rm as}\quad s\to0^+.
\end{eqnarray*}
The proof ends. \qquad$\Box$

\smallskip

From Lemma \ref{lm 3.1}, it is informed that the limit of $\mathbb{G}_\alpha[\Gamma_s]$ as $s\to0^+$ can't be used
as a barrier function to control the sequence $\{u_s\}$. So we have to
find new upper   bound for sequence $\{u_s\}$.

\begin{proposition}\label{pr 4.1}
 Let $\Gamma_0$ be defined in (\ref{Gamma}) and
\begin{equation}\label{5.4}
p>1+\frac{2\alpha}N,
\end{equation} then problem
\begin{equation}\label{5.5}
 \arraycolsep=1pt
\begin{array}{lll}
 (-\Delta)^{\alpha} u+u^p=\Gamma_0\ \ \ \ & {\rm in}\ \
B_1(e_N),\\[2mm]\phantom{(-\Delta)^{\alpha}+u^p}
u=0& {\rm in} \ \ B_1^c(e_N)\setminus\{0\}.
\end{array}
\end{equation}
admits a minimum positive solution $u_0$, that is, $u_0\le u$ for any nonnegative solution $u$ of (\ref{5.5}).

Moreover,
\begin{equation}\label{5.6}
\lim_{x\in B_1(e_N),\ x\to\partial B_1(e_N)\setminus\{0\}}u_0(x)=0
\end{equation} and
\begin{equation}\label{5.7}
\frac1{c_{12}}t^{-\frac{N+2\alpha}p}\le
u_0(te_N)\le c_{12}  t^{-\frac{N+2\alpha}p},\quad t\in(0,1),
\end{equation}
where  $c_{12}>1$ is independent of $\alpha$.
\end{proposition}
{\bf Proof. }
{\it The existence of solution to (\ref{5.5}).}   It implies by Lemma \ref{lm 3.1} that the mapping
$s\mapsto u_s$ is decreasing in $B_1(e_N)$, where $u_s\chi_{B_1(e_N)}$ is the solution of (\ref{eq 2.4}).
So what we have to do is just to find a super solution $U$ of (\ref{5.5})
such that $u_0\leq U$ in $ B_1(e_N)$. To this end, we consider the
radial function
\begin{equation}\label{02-10}
\Phi_\sigma(x)=\frac1{|x|^{\sigma}},\quad\forall x\in \R^N\setminus\{0\},
\end{equation}
where $\sigma\in[0,\ N)$.
By scaling property of $\Phi_\sigma$, (also see \cite{FQ2}) we know  that
\begin{equation}\label{5.8}
(-\Delta)^\alpha \Phi_\sigma(x)=\frac{c(\sigma,\alpha) }{|x|^{\sigma+2\alpha}},
\end{equation}
where  $c(\sigma,\alpha)\in\R$. Now we choose
 $\sigma=\sigma_0=\frac{N+2\alpha}p$, then $\sigma_0\in(0,\ N)$ if
$p>1+\frac{2\alpha}N$.   Therefore, there exist some $k>1$ dependent of $|c(\sigma_0,\alpha)|$ and $c_{N,\alpha}$ but independent of $n$ such that
\begin{equation}\label{17-09-1}
U(x)= k  \Phi_{\sigma_0}(x)
\end{equation}
is a super solution  of (\ref{5.5}). Thus, $U\in L^1( B_1(e_N))$ and
by Theorem \ref{comparison}, we have that
\begin{equation}\label{5.9}
0\leq u_s\leq U\ \ {\rm for\ any}\ \ s\in(0,1).
\end{equation}
For any $x\in\R^N\setminus\{0\}$, $
u_0(x):=\lim_{s\to0^+}u_s(x)\leq U(x)<+\infty$.  Following the same argument of Proposition \ref{pr 2.1}, we can prove that $ u_0$ is a classical  solution of (\ref{5.5}).
Furthermore, $u_0$ is the minimum solution of (\ref{5.5}).

\smallskip

{\it  Proof of (\ref{5.6}).} Let $\bar
x\in\partial B_1(e_N)\setminus\{0\}$, $K_1=\partial B_1(e_N)\cap \overline{B_{|\bar x|/8}(\bar x)}$ and $K_2=\partial B_1(e_N)\cap  B_{|\bar
x|/2}^c(\bar x)$.   Let $\mathcal{O}$ be an open and
$C^2$ set such that
$$B_1(e_N)\cap  B_{|\bar x|/4}(\bar x)\subset \mathcal{O}\subset B_1(e_N)\cap B_{|\bar x|/2}(\bar x).$$
Then we see that
$$K_1\subset \partial \mathcal{O}\quad{\rm and}\quad \partial \mathcal{O}\cap K_2=\emptyset.$$

 We would like to find
a super solution of (\ref{5.5}) in $\mathcal{O}$ with vanishing
boundary value in $K_1$ for any $n$.  Denote
$$V_\lambda =U\eta+\lambda  V_{\mathcal{O}},$$
where $U$ is defined (\ref{17-09-1}), $\eta$ is a $C^2$ function such that
$$
\eta(x)=\left\{ \arraycolsep=1pt
\begin{array}{lll}
 1\ \ \ \ &{\rm if}\quad
x\in B^c_{|\bar x|/4}(\bar x),\\[2mm]
0\ &{\rm if}\quad x\in B_{|\bar x|/8}(\bar x).
\end{array}
\right.
$$
  and
$ V_{\mathcal{O}}$ is the solution of
$$
 \arraycolsep=1pt
\begin{array}{lll}
 (-\Delta)^{\alpha}  u = 1\quad   &{\rm in}\quad \mathcal{O},\\[2mm]
 \phantom{(-\Delta)^{\alpha}}
 u=0\quad   &{\rm in}\quad \mathcal{O}^c.
\end{array}
$$
Since $U(1-\eta)$ is $C^2$ in $\R^N$, then there exists $c_{13}>0$ such that
\begin{eqnarray*}
|(-\Delta)^{\alpha} U(1-\eta)|\le c_{13}\quad {\rm in}\quad \bar\mathcal{O}.
\end{eqnarray*}
Thus, there exists $c_{14}>0$ such that
$|(-\Delta)^{\alpha} U\eta|\le c_{14}\quad {\rm in}\quad \bar\mathcal{O}.$
Choosing $\lambda_0>0$ suitable,  we have that
for $\lambda\ge \lambda_0$,
$$(-\Delta)^{\alpha}V_\lambda + V_\lambda^p \ge \frac{c_{N,\alpha}}{|x|^{N+2\alpha}}\ \ {\rm in}\ \ \mathcal{O}.$$
Moreover, since $\eta=1$ in $\R^N\setminus \mathcal{O}$, then $V_\lambda\ge U\ge u_s$ in $\R^N\setminus \mathcal{O}$.
By the Comparison Principle, we have that for any $s\in(0,1)$,
\begin{equation}\label{17-09-2}
u_s\le V_\lambda\quad{\rm in}\quad B_1(e_N).
\end{equation}
which implies (\ref{5.6}).

{\it Proof of (\ref{5.7}).}  For $t\in(0,1)$, denote that
$$V_t(x)=c_{N,\alpha}^{\frac1{p}}t^{-\frac{N+2\alpha}p} V_B\left(\frac{x-te_N}{t}\right), \quad x\in\R^N, $$
 where  $V_B$ is the solution of
\begin{equation}\label{eq 3.1}
 \arraycolsep=1pt
\begin{array}{lll}
 (-\Delta)^{\alpha}  u= 1\quad&{\rm in}\quad B_1(0),\\[2mm] \phantom{ (-\Delta)^{\alpha} }
 u=0\quad&{\rm in}\quad \R^N\setminus B_1(0).
\end{array}
\end{equation}

It deduces by (\ref{5.4}) that $x\in B_{\frac t4}(te_N)$,
$$(-\Delta)^{\alpha} V_t(x)=c_{N,\alpha}^{\frac1{p}}t^{-\frac{N+2\alpha}p-2\alpha}\le c_{N,\alpha}^{\frac1{p}}t^{-N-2\alpha} $$
and $$ \frac1{|x+te_N|^{N+2\alpha}}\ge \frac {c_{15}}{t^{N+2\alpha}},$$
where $c_{15}>0$ independent of $t$.
Then there exists some constant $\nu\in(0,1)$ such that
$$(-\Delta)^{\alpha}(\nu V_t)+(\nu V_t)^p\leq \frac{c_{N,\alpha}}{t^{N+2\alpha}}\leq (-\Delta)^{\alpha}u_0+u_0^p\quad {\rm in}\quad  B_{\frac t4}(te_N)$$
and $\nu V_t=0\leq u_0$ in $B_{\frac t4}(te_N)^c$. By applying the Comparison Principle, we have that
$$\nu V_t\leq u_0\quad {\rm in}\quad \R^N,$$
 which implies
\begin{equation}
 u_0(x)\ge \nu c_{N,\alpha}^{\frac1{p}} t^{-\frac{N+2\alpha}p}\min_{B_{\frac 12}(0)}V_B\ge
c_{16}c_{N,\alpha}^{\frac1{p}} |x|^{-\frac{N+2\alpha}p},\quad\
\forall x\in B_{\frac t8}(te_N),
\end{equation}
for some constants $c_{16}>0$ independent of $t$.
We complete the proof.\qquad$\Box$

\smallskip
 Now we are in the position to prove Theorem \ref{teo 2}.
\smallskip

\noindent{\bf Proof of Theorem \ref{teo 2}.}
 By Proposition \ref{pr 4.1},  there exists  minimum solution $u_0$ of (\ref{5.5}). Using the maximum principle argument,
 we have that for $s\in (0,1)$,
 $$u_s\le u_0\quad{\rm in}\quad  B_1(e_N).$$

 Since $$u_0\in L^1(B_1(e_N)),\quad u_0^p\in L^1(B_1(e_N),\rho^2 dx),$$
  where $\rho(x)={\rm dist}(x,\partial B_1(e_N))$,
   then there exists $\tilde u_0\le u_0$ such that
$$u_s\to \tilde u_0\quad {\rm in}\ \ B_1(e_N)\quad {\rm and} \quad {\rm in}\ \ L^1(B_1(e_N))$$
 and
 $$u_s^p\to \tilde u_0^p\quad {\rm in}\ \ B_1(e_N)\quad {\rm and} \quad {\rm in}\ \  L^1(B_1(e_N),\rho^2 dx).$$

Passing the limit in identity of (\ref{2.5}) as $s\to0^+$, we have that
\begin{equation}\label{17-09-0}
\int_{B_1(e_N)}[\tilde u_0(-\Delta)^\alpha \xi+\tilde u_0^p\xi]dx= \int_{B_1(e_N)} \xi \Gamma_0   dx,\qquad\forall \xi\in C_0^\infty(B_1(e_N)).
\end{equation}

 By the same argument, we have that $\tilde u_0$ is $C^2$ locally in $B_1(e_N)$ and $\tilde u_0\chi_{B_1(e_N)}$ is the minimum solution of
 (\ref{5.5}).
  Thus, we have
 $$\tilde u_0=u_0\quad{\rm in}\quad B_1(e_N). $$

{\it Uniqueness.} Let $v_0$ be a weak solution of (\ref{eq 2.2}) and then $\varphi_0:=u_0-v_0$ is a weak solution of
$$
\arraycolsep=1pt
\begin{array}{lll}
(-\Delta)^\alpha \varphi_0 +u_0^p- v_0^p=0\quad &{\rm in}\quad B_1(e_N),
\\[2mm]\phantom{(-\Delta)^\alpha   +u_{s}^p- v_s^p}
\varphi_0=0\quad &{\rm in}\quad \R^N\setminus B_1(e_N),
\end{array}
$$
that is,
\begin{equation}\label{25-09-0}
\int_{ B_1(e_N)}\varphi_0(-\Delta)^\alpha\xi dx+\int_{ B_1(e_N)}(u_0^p- v_0^p)\xi dx=0,\qquad \forall \xi\in C_0^\infty(B_1(e_N)).
\end{equation}
In this definition of weak solution, we can not apply Kato's inequality \cite[Proposition 2.4]{CV1} directly due to the stronger test functions space.
So we have to improve the regularity of $\varphi_0$. In fact, for any $\xi\in\mathbb{X}_{\alpha,B_1(e_N)}$, let $\{\xi_k\}$ be a sequence of nonnegative functions in $C_0^\infty(B_1(e_N))$,
such that
$$\xi_k\to \xi\quad {\rm and}\quad |(-\Delta)^\alpha\xi_k|\le 2|(-\Delta)^\alpha\xi|\quad {\rm in}\quad B_1(e_N). $$
Passing the limit of (\ref{25-09-0}) with $\xi_k$ as $k\to\infty$, we have that
(\ref{25-09-0}) holds for any $\xi\in \mathbb{X}_{\alpha,B_1(e_N)}.$
If we choose a sequence $\{\tilde \xi_k\}\subset \mathbb{X}_{\alpha,B_1(e_N)}$ which converges
to ${\rm sign}(u_0^p- v_0^p)\cdot \mathbb{G}_{\alpha,B_1(e_N)}[1]$ and then there exists $c_6>0$ such that
$$\int_{ B_1(e_N)}|u_0^p- v_0^p|\mathbb{G}_{\alpha,B_1(e_N)}[1] dx\le c_6\|\varphi_0\|_{L^1(B_1(e_N))},$$
therefore, $|u_0^p- v_0^p|\in L^1(B_1(e_N),\rho^\alpha dx)$.

Now we follows the Kato's inequality \cite[Proposition 2.4]{CV1} to obtain that
$$\int_{B_1(e_N)}|\varphi_0|(-\Delta)^\alpha\xi+\int_{B_1(e_N)}[ u_0^p-v_0^p]{\rm sign}(u_0-v_0)\xi dx=0,\qquad\forall \xi\in\mathbb{X}_{\alpha,B_1(e_N)}. $$
Taking $\xi=\mathbb{G}_\alpha[1]$, we have that
$$\int_{B_1(e_N)}[u_0^p- v_0^p]{\rm sign}(u_0-v_0)\xi dx\ge 0\quad{\rm and}\quad \int_{B_1(e_N)}|\varphi_0|dx=0,$$
then $\varphi_0=0$ a.e. in $B_1(e_N)$. Then the uniqueness is proved. \qquad$\Box$

 \setcounter{equation}{0}
\section{ Nonexistence}

In order to prove  the nonexistence of weak solution to (\ref{1.1}) for $p\in(0,1+\frac{2\alpha}N]$, we will prove that
the solution $u_s$ of problem (\ref{eq 2.1}) blows up in $B_1(e_N)$ as $s\to0^+$.
To this end, we first define a cone  by
\begin{equation}\label{cone}
\mathcal{C}=\{x\in\R^N:\  \exists t\in(0,1)\ {\rm s.t.}\  |x-te_N|<\frac t8 \}.
\end{equation}
We observe that $\mathcal{C}\subset B_1(e_N)$.
\begin{lemma}\label{lm 4.1}
Assume that $p>0$, $s\in(0,1)$ and  $u_s$ is the unique solution of problem (\ref{eq 2.1}).
Then there exists $c_{17}>0$ such that
$$\lim_{s\to0^+}u_s(x)\ge c_{17}|x|^{-N},\quad \forall x\in \mathcal{C}.$$

\end{lemma}
{\bf Proof.} By Proposition \ref{pr 2.1}, we have that
 $\tilde u_s:=u_s\chi_{B_1(e_N)}$ is the unique classical solution of (\ref{eq 2.4}).
  For given $t\in(0,1)$, we recall
$$V_t(x)= t^{-N} V_B\left(\frac{x-te_N}{t}\right), $$
 where  $V_B$ is the solution of (\ref{eq 3.1}).

We see that for  $x\in B_{\frac t4}(te_N)$,
$$(-\Delta)^{\alpha} V_t(x)= t^{-N-2\alpha}$$
and   choose $s\le \frac t4$, then
$$|x+se_N|\le |x|+s\le \frac32 t,$$
which implies that
$$ \frac1{|x+se_N|^{N+2\alpha}}\ge \frac{c_{18}}{t^{N+2\alpha}},$$
where $c_{18}>0$ independent of $t$.
Then there exists some constant $\nu>0$ such that
$$(-\Delta)^{\alpha}(\nu V_t)+(\nu V_t)^p\leq \frac {c_{N,\alpha}}{t^{N+2\alpha}}\leq (-\Delta)^{\alpha}\tilde u_s+\tilde u_s^p\quad {\rm in}\quad  B_{\frac t4}(te_N)$$
and $\nu V_t=0\leq \tilde u_s$ in $B_{\frac t4}(te_N)^c$, by applying the Comparison Principle, we have that
$$V_t\leq \tilde u_s\quad {\rm in}\quad \R^N,$$
 which implies
\begin{equation}\label{5.10}
 u_s(x)\ge \nu   t^{-N}\min_{B_{\frac 12}(0)}V_B\ge
c_{19}  |x|^{-N},\quad \forall
x\in B_{\frac t8}(te_N),
\end{equation}
for some constants $c_{19}>0$ independent of $t$ and $s$.
Combining the increasing monotonicity, we have that
$$\lim_{s\to0^+}u_s(x)\ge c_{19}|x|^{-N},\quad \forall x\in \mathcal{C}.$$

\begin{lemma}\label{lm 4.2}
Assume that $p\ge 0$, $s\in(0,1)$ and  $\tilde u_s:=u_s\chi_{B_1(e_N)}$ is the minimal solution of (\ref{eq 2.4}).
 Then \\$(i)$\ \ $\tilde u_s(x',x_N)$ is symmetric with respect to $x'$ and decreasing with to $r':=|x'|$ for any $x_N\in (0,2)$;\\
 $(ii)$  $\tilde u_s(x',x_N)$ is decreasing in $x_N$ for $x_N\in(1,1+\sqrt{1-|x'|^2})$.

\end{lemma}
{\bf Proof.}   By applying the same procedure of step 1 and step 2 in  proof of Theorem 1.1 in \cite{FW},   we have that
 $$u_s(x', x_N)= u_s(|x'|, x_N)\quad{\rm for}\ \ x=(x', x_N)\in B_1(e_N)$$
 and $u_s(r, x_N)$ is decreasing with   $r=|x'|$.
Using the same argument from the other side, we conclude that $u_s(x',x_N)\le u_s(x',2-x_N)$ for $ x_N\in(1,2)$,
that is, $u_s(x',x_N)$ is decreasing with $x_N\in(1,2)$.

 \hfill$\Box$

\medskip

\noindent{\bf Proof of Theorem \ref{teo 3}.}
For $p=0$, the solution $u_s$ of (\ref{eq 2.1}) satisfies
$$u_s=\mathbb{G}_\alpha[\Gamma_s]-\mathbb{G}_\alpha[1].$$
By Lemma \ref{lm 3.1}, we know that $\mathbb{G}_\alpha[\cdot]$ blows every where in $B_1(e_N)$,
but $\mathbb{G}_\alpha[1]$ is bounded uniformly.

For $0<p\le 1+\frac{2\alpha}N$, our proof is divided into two steps.

{\it Step 1. we prove that
$$\lim_{s\to0^+}u_s(x)=+\infty,\quad x\in \mathcal{B}:=\{(x',x_N)\in B_1(e_N):\ x_N>1\}.$$}
If there is $\bar x=(\bar x_0, t_0)\in \mathcal{B}$ such that $\lim_{s\to0^+}u_s(\bar x)<+\infty$, then by Lemma \ref{lm 4.2}, we have that
$$\lim_{s\to0^+}u_s(x)\le \lim_{s\to0^+}u_s(\bar x)\quad {\rm in }\quad B_+:=\{x=(x',x_N)\in B_1(e_N):\ |x'|>|\bar x'|,\ x_N>t_0\}.$$
Choose a nonnegative function $\xi_0\in C_0^\infty$ with support in $B_+$, then for any $x\in B_1(e_N)\setminus B_+$,
\begin{eqnarray*}
(-\Delta)^\alpha \xi_0(x) &=& -c_{N,\alpha}{\rm P.V.}\int_{\R^N}\frac{\xi_0(y)-\xi_0(x)}{|x-y|^{N+2\alpha}}dy \\
   &=& -c_{N,\alpha} \int_{ B_+}\frac{\xi_0(y)}{|x-y|^{N+2\alpha}}dy
  \\&\le & -c_{N,\alpha} |x-2e_N|^{-N-2\alpha}\int_{ B_+}\xi_0(y)dy,
\end{eqnarray*}
therefore, there exists $c_{20}>0$ such that
$$(-\Delta)^\alpha \xi_0(x)\le -c_{20},\quad\forall x\in B_1(e_N)\setminus B_+.$$
Moreover, there exists $c_{21}>0$ such that
$$|(-\Delta)^\alpha \xi_0(x)|\le c_{21},\quad\forall x\in  B_+.$$
Since $\mathcal{C}\subset (B_1(e_N)\setminus B_+)$, then by Lemma \ref{lm 4.1} we have that
\begin{eqnarray*}
 \int_{B_1(e_N)}u_s(-\Delta)^\alpha\xi_0dx &\le & \int_{\mathcal{C}}u_s(-\Delta)^\alpha\xi_0dx+\int_{B_1(e_N)\setminus (B_+\cup \mathcal{C})} u_s(-\Delta)^\alpha\xi_0dx\\&&+\int_{B_+}u_s|(-\Delta)^\alpha\xi_0|dx
 \\ &\le &-c_{20}  \int_{  \mathcal{C}}  u_s dx+c_{21}|B_+|
 \\&\to&-\infty\quad{\rm as}\quad s\to0^+.
\end{eqnarray*}
While
\begin{eqnarray*}
  |\int_{B_1(e_N)}u_s^p\xi_0dx| = \int_{B_+}u_s^p\xi_0dx    \le u_s^p(\bar x)\max_{\R^N}\xi_0 |B_+|
\end{eqnarray*}
 and
\begin{eqnarray*}
  \int_{B_1(e_N)}\frac{\xi_0(x)}{|x+se_N|^{N+2\alpha}}dx = \int_{B_+}\frac{\xi_0(x)}{|x+se_N|^{N+2\alpha}}dx  \le  |B_+|\max_{\R^N}\xi_0,
\end{eqnarray*}
where $|x+se_N|\ge 1$ for $x\in B_+$.
Then taking $s>0$ small enough, we obtain a contradiction with
 the identity
$$\int_{B_1(e_N)}[u_s(-\Delta)^\alpha \xi+u_s^p\xi]dx=c_{N,\alpha}\int_{B_1(e_N)} \frac{\xi(x)}{|x+se_N|^{N+2\alpha}} dx,\qquad\forall \xi\in C_0^\infty(B_1(e_N)).$$
Therefore,
$$\lim_{s\to0^+}u_s(x)=+\infty,\quad \forall x\in \mathcal{B}.$$

 {\it Step 2: We claim that $\lim_{s\to0^+}u_s(x)=\infty$, $x\in B_1(e_N)$}. By the  fact of $\lim_{s\to0^+}u_s(x)=\infty$, $x\in \mathcal{B}$,
then letting $\tilde x=(\frac12,0\cdots, \frac32)$,  for any $n>1$,
there exists $s_n>0$ such that $s_n\to0$ as $n\to+\infty$ and
$$u_{s_n}(\tilde x)\ge n,$$
then applying Lemma \ref{lm 4.2}, we have that
$$u_{s_n}\ge n\quad{\rm in}\quad \mathcal{B}_0=\left\{x=(x',x_N)\in\mathcal{B}:\ |x'|\le \frac12,\ 1\le x_N\le \frac32\right\}.$$

For any $x_0\in B_1(e_N)\setminus \mathcal{B}$, there exists ${r_1}>0$ such that
$\bar B_{r_1}(x_0)\subset B_1(e_N)\setminus \mathcal{B}$.
 We denote by $\phi_{n}$ the solution of
\begin{equation}
\arraycolsep=1pt
\begin{array}{lll}
 (-\Delta)^\alpha  u+u^p=0 \quad & {\rm in}\quad  B_{r_1}(x_0),\\[2mm]
 \phantom{  (-\Delta)^\alpha  +u^p}
u=0  \quad & {\rm in}\quad B_{r_1}^c(x_0)\setminus \mathcal{B}_0,\\[2mm]
\phantom{ (-\Delta)^\alpha  +u^p}
u=n  \quad & {\rm in}\quad \mathcal{B}_0.
\end{array}
\end{equation}
Then by Theorem \ref{comparison}, we have that
\begin{equation}\label{4.1.3}
 u_{s_n}\ge \phi_n\quad {\rm in}\quad B_1(e_N).
\end{equation}
Let  $\varphi_n=\phi_n-n \chi_{\mathcal{B}_0},$
then $\varphi_n=\phi_n$ in $B_{r_1}(x_0)$ and
\begin{eqnarray*}
      (-\Delta)^\alpha \varphi_n(x)+\varphi_n^p(x) &=&  (-\Delta)^\alpha \phi_n(x)-  n (-\Delta)^\alpha \chi_{\mathcal{B}_0}(x) +\phi_n^p(x) \\
        &=& n \int_{\mathcal{B}_0}\frac{dy}{|y-x|^{N+2\alpha}},\qquad \forall x\in B_{r_1}(x_0),
     \end{eqnarray*}
that is,  $\varphi_n$  is a solution of
\begin{equation}\label{4.1.2}
\arraycolsep=1pt
\begin{array}{lll}
\displaystyle (-\Delta)^\alpha  u+u^p= n \int_{\mathcal{B}_0}\frac{dy}{|y-x|^{N+2\alpha}} \quad & {\rm in}\quad  B_{r_1}(x_0),\\[2mm]
 \phantom{  (-\Delta)^\alpha  +u^{p,}}
u=0  \quad & {\rm in}\quad B_{r_1}^c(x_0).
\end{array}
\end{equation}
By direct computation,
$$\frac1{c_{23}}\le \int_{\mathcal{B}_0}\frac{dy}{|y-x|^{N+2\alpha}}\le c_{23},\quad\forall x\in B_{r_1}(x_0),$$
for some $c_{23}>1$.

Let $\eta_1$ be the solution of
$$
\arraycolsep=1pt
\begin{array}{lll}
 (-\Delta)^\alpha  u=1 \quad & {\rm in}\quad  B_{r_1}(x_0),\\[2mm]
 \phantom{  (-\Delta)^\alpha  }
u=0  \quad & {\rm in}\quad  B^c_{r_1}(x_0)
\end{array}
$$
and then
 $(\frac{ n}{2c_{23}})^{\frac1p} \max\eta_1\cdot\eta_1$ is sub solution of (\ref{4.1.2}) for $ n$ large enough. Then it infers by Theorem \ref{comparison} that
 $$\varphi_n\ge (\frac{ n}{2c_{23}})^{\frac1p} \max\eta_1\cdot\eta_1,\qquad \forall  x\in B_{r_1}(x_0),$$
 which implies that
 $$\phi_n\ge (\frac{ n }{2c_{23}})^{\frac1p}\max\eta_1\cdot \eta_1,\qquad \forall  x\in B_{r_1}(x_0).$$
 Then by (\ref{4.1.3}), $$\lim_{s\to0^+}u_{s_n}(x_0)\ge \lim_{n\to\infty}\phi_n(x_0)=\infty.$$
Since $x_0$ is arbitrary in $B_1(e_N)\setminus \mathcal{B}$, it implies that $\lim_{n\to\infty^+}u_{s_n}(x)=\infty$ in $B_1(e_N)$.\qquad$\Box$

\setcounter{equation}{0}
\section{The solutions vanishes as $\alpha\to1^-$ }

\subsection{The case of $p\ge \frac{N+2}{N-2}$ }
\begin{lemma}\label{lm 2.0}
$(i)$ Let $N\ge 2$ and  $c_{N,\alpha}$ define in (\ref{1.2}). Then
$$\lim_{\alpha\to1^-}\frac{c_{N,\alpha}}{1-\alpha}=\frac{4N}{|S^{N-1}|},$$
where $|S^{N-1}|$ denotes the $(N-1)-$dimensional measure of the unit sphere $S^{N-1}$.

$(ii)$ Let $N\ge 2$, then for any $f\in C_0^\infty(\R^N)$,
$$\lim_{\alpha\to1^-}(-\Delta)^\alpha f=-\Delta f. $$

\end{lemma}
{\bf Proof.} The proofs of $(i)$ and $(ii)$ see Corollary 4.2 and Proposition 4.4 in \cite{NPV} respectively.\qquad$\Box$
\medskip

\begin{proposition}\label{pr 2.00}
Let $\Phi_\sigma$ be defined by (\ref{02-10}) with $\sigma\in(0,\ N)$ and $$c(\sigma,\alpha)=-\frac12\int_{\R^N}\frac{|z+e_N|^{-\sigma}+|z-e_N|^{-\sigma}-2}{|z|^{N+2\alpha}}dz.$$
Then
\begin{equation}\label{22-09-0}
(-\Delta)^\alpha \Phi_\sigma(x)=\frac{c(\sigma,\alpha)}{|x|^{\sigma+2\alpha}},\quad\forall x\in \R^N\setminus\{0\}.
\end{equation}
Moveover,
$$\lim_{\alpha\to1^-}c(\sigma,\alpha)=(N-2-\sigma)\sigma$$
and
 \begin{equation}\label{22-09-1}
\lim_{\alpha\to1^-}(-\Delta)^\alpha \Phi_\sigma(x)=\frac{(N-2-\sigma)\sigma}{|x|^{\sigma+2}},\quad\forall x\in \R^N\setminus\{0\}.
\end{equation}

\end{proposition}
{\bf Proof.} By direct computation, we have that
\begin{equation}\label{22-09-2}
(-\Delta)^\alpha \Phi_\sigma(x)=\frac{c(\sigma,\alpha)}{|x|^{\sigma+2\alpha}},\quad\forall x\in \R^N\setminus\{0\},
\end{equation}
where
$$c(\sigma,\alpha)=-\frac12\int_{\R^N}\frac{\frac1{|z+e_N|^{\sigma}}+\frac1{|z-e_N|^{\sigma}}-2}{|z|^{N+2\alpha}}dz.$$
On the other hand, it has been proved in \cite{CS1} that
$$c(\sigma,\alpha)=0\quad{\rm if}\quad  \sigma=N-2\alpha.$$

Now for any $R>1$, let $\eta_R:\R^N\to [0,1]$ be a nonnegative function such that $\eta_R=1$ in $B_R(0)\setminus B_{\frac1R}(0)$
and $\eta_R=0$ in $B_{2R}^c(0)\cup B_{\frac1{2R}}(0)$, then $\eta_R\Phi_\sigma\in C^\infty_0(\R^N)$. By Lemma \ref{lm 2.0} $(ii)$,
we have
$$\lim_{\alpha\to1^-}(-\Delta)^\alpha (\eta_R\Phi_\sigma)(e_N)=(N-2-\sigma)\sigma $$
and
\begin{equation}\label{02-10-0}
(-\Delta)^\alpha (\eta_R\Phi_\sigma)(e_N)=(-\Delta)^\alpha \Phi_\sigma(e_N)-c_{N,\alpha}\int_{\R^N}\frac{(1-\eta_R(z))|z|^{-\sigma}}{|z-e_N|^{N+2\alpha}}dz.
\end{equation}
A straightforward computation implies that
\begin{eqnarray*}
  0 &\le & \int_{\R^N}\frac{(1-\eta_R(z))|z|^{-\sigma}}{|z-e_N|^{N+2\alpha}}dz \\
   &\le &  \int_{B_{\frac1R}(0)\cup B_R^c(0)}\frac{|z|^{-\sigma}}{|z-e_N|^{N+2\alpha}}dz
  \\ &\le & c_{24}\int_{B_{\frac1R}(0)\cup B_R^c(0)}\frac{|z|^{-\sigma}}{|z-e_N|^{N}}dz,
\end{eqnarray*}
where $c_{24}>0$ is independent of $\alpha$.
Thus,
$$\lim_{\alpha\to1^-}c_{N,\alpha}\int_{\R^N}\frac{(1-\eta_R(z))|z|^{-\sigma}}{|z-e_N|^{N+2\alpha}}dz=0.$$
Therefore, passing to the limit of (\ref{02-10-0}) as $\alpha\to1^-$,  we conclude that $\lim_{\alpha\to1^-}c(\sigma,\alpha)=(N-2-\sigma)\sigma$
and combining (\ref{22-09-2}),  (\ref{22-09-1}) holds.\qquad$\Box$

\begin{proposition}\label{pr 5.1}
 Let $p\ge \frac{N+2}{N-2}$, $\sigma_p=\frac{N+2}p$ and
$u_{\alpha,p}$ be the unique weak solution of  (\ref{1.1}).
Then
\begin{equation}\label{5.1}
u_{\alpha,p}\le (4^{1-\alpha}c_{N,\alpha})^{\frac1p}\Phi_{\sigma_p}\quad {\rm in }\ B_1(e_N).
\end{equation}
\end{proposition}
{\bf Proof. }
It follows by (\ref{22-09-0}) that
\begin{equation}\label{5.2}
(-\Delta)^\alpha \Phi_{\sigma_p}(x)=\frac{c(\sigma_p,\alpha)}{|x|^{\sigma_p+2\alpha}},\quad \forall x\in\R^N\setminus\{0\},
\end{equation}
where $c(\sigma_p,\alpha)\to N-2-\sigma_p$ as $\alpha\to 1^-$.
For $k>1$, denote
$$
U_k= k  \Phi_{\sigma_p}.
$$
When $p> \frac{N+2}{N-2}$, we have $N-2-\sigma_p>0$, then there exists $\alpha_p\in(0,1)$ such that $c(\sigma_p,\alpha)\ge 0$ for $\alpha\in(\alpha_p,1)$.
When $p=\frac{N+2}{N-2}$, we have $\sigma_p=N-2<N-2\alpha$. Using \cite{FQ2}, $c(\cdot,\alpha)$ is $C^2$ and convex in $[0,N)$ and
$$\lim_{\sigma\to N^-}c(\sigma,\alpha)=-\infty$$
and for $\sigma\in(0,N-2\alpha)$,
$$ c(\sigma,\alpha)\ge c(N-2\alpha,\alpha)=0.$$
Thus, $c(\sigma_p,\alpha)\ge0$ by the fact of $\sigma_p<N-2\alpha$.
Choosing $k=(4^{1-\alpha}c_{N,\alpha})^{\frac1p}$, then for any $x\in B_1(e_N) $,
\begin{eqnarray*}
(-\Delta)^\alpha U_k+U_k^p-\frac{c_{N,\alpha}}{|x|^{N+2\alpha}} &=& \frac{kc(\sigma_p,\alpha)}{|x|^{\sigma_p+2\alpha}}+\frac{k^p}{|x|^{N+2}} -\frac{c_{N,\alpha}}{|x|^{N+2\alpha}}\\
   &\ge & (k^p -4^{1-\alpha}c_{N,\alpha})\frac{1}{|x|^{N+2}}=0,
\end{eqnarray*}
where we used that $\frac{1}{|x|^{N+2\alpha}}\le \frac{4^{1-\alpha}}{|x|^{N+2}}$ for $x\in B_1(e_N)$.
Therefore, for $\alpha\in(\alpha_p,1)$ we have that
$$u_{\alpha,p}(x)\le(4^{1-\alpha}c_{N,\alpha})^{\frac1p} \Phi_{\sigma_p}(x),\qquad \forall x\in\ B_1(e_N).\qquad\Box$$

\medskip
\noindent{\bf Proof of Theorem \ref{teo 4}.} It follows by Proposition \ref{pr 5.1} that
$$u_{\alpha,p}\le (4^{1-\alpha}c_{N,\alpha})^{\frac1p}\Phi_{\sigma_p}\quad {\rm in }\ B_1(e_N).$$
On the other hand,  we deduce from Lemma \ref{lm 2.0} that
$$c_{N,\alpha}\le c_{25}(1-\alpha),$$
where $c_{25}>0$ is independent of $\alpha$.
Therefore, we have that
$$0\le \lim_{\alpha\to1^-}u_{\alpha,p}(x)\le c_{26}\lim_{\alpha\to1^-} (1-\alpha)^{\frac1p}\Phi_{\sigma_p}(x)=0,\quad \forall x\in B_1(e_N),$$
where $c_{26}>0$ is independent of $\alpha$. This ends the proof.\qquad $\Box$

\subsection{The case $p\in [\frac{N+1}{N-1}, \frac{N+2}{N-2})$ }

We know that the unique weak solution $u_{\alpha,p}$ of (\ref{1.1}) satisfies
the identity
\begin{equation}\label{15-10-0}
\int_{B_1(e_N)} \left[u_{\alpha,p}(x)(-\Delta)^\alpha \xi(x) +u_{\alpha,p}^p(x)\xi(x)\right]dx =
\int_{B_1(e_N)}\xi(x)\Gamma_{0,\alpha}(x) dx,\quad \forall \xi\in C^\infty_0(B_1(e_N)).
\end{equation}
From Lemma \ref{lm 2.0} part $(ii)$, it infers that for any
$\xi\in C^\infty_0(B_1(e_N))$,
$$(-\Delta)^\alpha \xi \to-\Delta\xi\quad{\rm uniformly\ in}\quad  B_1(e_N) \quad {\rm as}\quad \alpha\to1^-.  $$

\begin{lemma}\label{lm 5.1}
 Let $\frac{N+1}{N-1}\le p< \frac{N+2}{N-2}$, $\sigma_p=\frac{N+2}p$ and
$u_{\alpha,p}$ be the unique weak solution of  (\ref{1.1}).
Then there exists $c_{27}>0$ independent of $\alpha$ such that
\begin{equation}\label{5.11}
0\le u_{\alpha,p}\le c_{27}\Phi_{\sigma_p}\quad {\rm in }\ B_1(e_N).
\end{equation}
\end{lemma}
{\bf Proof.} For $\frac{N+1}{N-1}\le p< \frac{N+2}{N-2}$, we have $\sigma_p<N$. It follows by (\ref{22-09-0}) that
$$
(-\Delta)^\alpha \Phi_{\sigma_p}(x)=\frac{c(\sigma_p,\alpha)}{|x|^{\sigma_p+2\alpha}},\quad \forall x\in\R^N\setminus\{0\},
$$
where $c(\sigma_p,\alpha)\to N-2-\sigma_p$ as $\alpha\to 1^-$. There exists $\alpha_p\in(0,1)$ such that
$$|c(\sigma_p,\alpha)|\le 2|N-2-\sigma_p|,\quad \alpha\in(\alpha_p,1).$$
For $k>1$, denote
$$
U_k= k  \Phi_{\sigma_p}.
$$

For $\alpha\in(\alpha_p,1)$ and any $x\in B_1(e_N) $, we deduce that
\begin{eqnarray*}
(-\Delta)^\alpha U_k+U_k^p-\frac{c_{N,\alpha}}{|x|^{N+2\alpha}} &\ge& -\frac{2k|N-2-\sigma_p|}{|x|^{\sigma_p+2\alpha}}+\frac{k^p}{|x|^{N+2}} -\frac{c_{N,\alpha}}{|x|^{N+2\alpha}}\\
   &\ge & \left(k^p -4^{1-\alpha}c_{N,\alpha}-2^{N+3-\sigma_p}k|N-2-\sigma_p|\right)\frac{1}{|x|^{N+2}},
\end{eqnarray*}
which implies that there exists $k_0>0$ independent of $\alpha$ such that
for $\alpha\in(\alpha_p,1)$ we have that
$$u_{\alpha,p}(x)\le k_0 \Phi_{\sigma_p}(x),\qquad \forall x\in\ B_1(e_N).\qquad\Box$$

\smallskip
 We have that
$\Phi_{\sigma_p}\in L^1(B_1(e_N))$ and
$$\norm{\Phi_{\sigma_p}}_{M^{p^*}(B_1(0),dx)}\le c_{28}, $$
where  $c_{28}>0$, $p^*=\frac{N}{\sigma_p}>1$ and $M^{p^*}(B_1(0))$ is the Marcinkiewicz space of
exponent $p^*$. By regularity results, we have that
$$\norm{u_{\alpha,p}}_{C^\beta_{\rm loc}(B_1(0))}\le c_{29},$$
where $c_{29}>0$ is independent of $\alpha$ by (\ref{5.11}).
Therefore, by compactness, we only have to prove that 0 is the only accumulation point of the sequence $\{u_{\alpha,p}\}_\alpha$.
  Let $u^*$ in $L^1(B_1(e_N))$ be accumulation point of the sequence $\{u_{\alpha,p}\}_\alpha$ and a subsequence, still denote $\{u_{\alpha,p}\}_\alpha$,
converge to $u^*$.

On the one hand, $\Gamma_{0,\alpha}$ converges to 0 uniformly in any compact subset of $B_1(e_N)$. Therefore, for $\xi\in C_0^\infty(B_1(0))$,
we have that
$$\lim_{\alpha\to1^-}\int_{B_1(e_N)}\xi(x)\Gamma_{0,\alpha}(x) dx=0.$$
Thus, we have that $u^*$ satisfies the identity
\begin{equation}\label{15-10-1}
\int_{B_1(e_N)} \left[u^*(x)(-\Delta) \xi(x) +(u^*)^p(x)\xi(x)\right]dx =0,\quad \forall \xi\in C^\infty_0(B_1(e_N)),
\end{equation}
On the other hand, by regularity result, we have that
$u^*\in C^2(B_1(e_N))\cap L^1(B_1(e_N))$ continuous up to boundary $\partial B_1(e_N)\setminus \{0\}$,
therefore,
$u^*$ is a nonnegative classical solution of
\begin{equation}\label{15-10-2}
 \arraycolsep=1pt
\begin{array}{lll}
 -\Delta  u+u^p=0\ \ \ \ & {\rm in}\ \
B_1(e_N),\\[2mm]\phantom{-\Delta +u^p}
u=0& {\rm on} \ \ \partial B_1(e_N)\setminus\{0\}.
\end{array}
\end{equation}

Then from Theorem 3.1 in \cite{GV}, we have that
$$u^*\in L^\infty(B_1(0))$$
and   Theorem 3.1 in \cite{GV},   $u^*$ is a classical solution of
\begin{equation}\label{15-10-3}
 \arraycolsep=1pt
\begin{array}{lll}
 -\Delta  u+u^p=0\ \ \ \ & {\rm in}\ \
B_1(e_N),\\[2mm]\phantom{-\Delta +u^p}
u=0& {\rm on} \ \ \partial B_1(e_N),
\end{array}
\end{equation}
where $p\in(1,\frac{N+1}{N-1})$.
By Strong Maximum Principle, we have that
$$u^*\equiv0\quad {\rm in}\quad B_1(e_N).$$

\noindent{\bf Acknowledgements:} The authors want to thank the referee for carefully reading the paper and making suggestions that
resulted in a great improvement in the clarity of the proofs. H. Chen is supported by National Natural Science Foundation
of China,  No:11401270
and
the Project-sponsored by SRF for ROCS, SEM.

\end{document}